\documentclass[11pt]{article}

\usepackage{amssymb,amsmath,amsthm}

\usepackage{graphicx}

\usepackage{color}
\definecolor{darkgreen}{rgb}{0,.4,0.2}
\definecolor{darkagenta}{rgb}{.5,0,.5}
\definecolor{darkred}{rgb}{0.85,0,0}%was0.85
\definecolor{darkblue}{rgb}{0,0,.6}
\definecolor{lightgray}{gray}{.95}
\definecolor{rgrey}{rgb}{.8,0.4,.4}  % faint
\definecolor{grey}{rgb}{.13,.13,.13}  % almost black

\newtheorem*{theorem*}{Theorem}

\advance \topmargin by -\headheight
\advance \topmargin by -\headsep
\evensidemargin \oddsidemargin
\begin{document}

\begin{centering}
{\Large \textbf{How to add two natural numbers in base phi}}

\bigskip

{\bf \large F.~Michel Dekking}

\bigskip

{ DIAM,  Delft University of Technology, Faculty EEMCS,\\ P.O.~Box 5031, 2600 GA Delft, The Netherlands.}

\medskip

{\footnotesize \it Email:  F.M.Dekking@TUDelft.nl}

\end{centering}

\medskip

\begin{abstract}
 \noindent In the base phi representation any natural number is written uniquely as a sum of powers of the golden mean with coefficients 0 and 1, where it is required that the product of two consecutive digits is always 0. In this self-contained paper we give a new, and short proof of the recursive structure of the base phi representations of the natural numbers.
\end{abstract}

\medskip


\bigskip

\section{Introduction}

  Base phi representations were introduced by George Bergman in 1957 (\cite{Bergman}). Let the golden mean be given by  $\varphi:=(1+\sqrt{5})/2$.\\
  Ignoring leading and trailing zeros, any  natural number $N$ can be written uniquely as
  $$N= \sum_{i=-\infty}^{\infty} d_i \varphi^i,\vspace*{-.0cm}$$
  with digits $d_i=0$ or 1, and where $d_id_{i+1} = 11$ is not allowed.
  As usual, we denote the base phi representation of $N$ as $\beta(N)$, and we  write these representations with a `decimal' point as
  $$\beta(N) = d_{L}d_{L-1}\dots d_1d_0\cdot d_{-1}d_{-2} \dots d_{R+1}d_R.$$
We give ourselves the freedom to write also non-admissible representations in this notation. For example, since $4 =2\times 2$ and $\beta(2)=10\cdot 01$, we will write $\beta(4)\doteq 20\cdot02$. Here the $\doteq$-sign indicates that we consider a non-admissible representation.

\medskip

  Our concern will be the recursive structure of the set of all numbers in their base phi representation. One can say that this structure was discovered in the series of papers\footnote{N.B.: these authors write the representations in reverse order}  \cite{Filliponi-Hart}, \cite{Hart98},  \cite{Hart-Sanchis-99}, and \cite{San-San}. \,
  A version of the recursive structure theorem is given in Proposition 3.1 and Proposition 3.2 in \cite{San-San}.
  Referring to these two propositions  the authors state: ``The full result is expressed in the following propositions, and was proved in Lemma 3.8 of \cite{Hart-Sanchis-99}". However, Lemma 3.8 in \cite{Hart-Sanchis-99} consists of eleven statements, all (except the rather trivial number (11)) about frequencies of occurrences of 1's and 0's. This means that, at least formally, there is no proof of the recursive structure theorem. We will fill this gap in the Section \ref{sec:3}.
  
  Finally, we mention that the recursive structure theorem plays a crucial role in the papers \cite{Dekk-phi-FQ} and \cite{Dekk-phi-sum}.

\section{Adding two base phi numbers}\label{sec:add}

We first mention that the natural number 2 has representation $\beta(2)=10\cdot01,$ since $\varphi + \varphi^{-2} =2.$ That this is correct, can be computed, using the equation  $\varphi^2=  \varphi+1.$ With some more work one finds that $\beta(4)=101\cdot01$.

A more convenient way to find the $\beta$-representations is to add $\beta(1)=1\cdot$\, repeatedly.\\
When we add two base phi numbers, then, in general, there is a carry both to the left and (two places) to the right:
$$\beta(5) = \beta(4+1)\doteq \beta(4)+\beta(1)= 101\cdot01 + 1\cdot \doteq 102\cdot01 \doteq 110\cdot02 = 1000\cdot1001.$$
Here we used twice that $2\varphi^n=\varphi^{n+1}+\varphi^{n-2}$ for all integers $n$, a direct consequence of $\beta(2)=10\cdot01.$\, 
Note that there is not only a {\it double carry}, but that we also have to get rid of the 11's, by replacing them with 100's.
This is allowed because of the equation $\varphi^{n+2}=  \varphi^{n+1}+\varphi^{n}.$ We call this operation a {\it golden mean shift}.

\medskip

\noindent  For the convenience of the reader we provide a list of the base phi representations of the first 24 natural numbers:

\medskip

 \begin{tabular}{|r|c|}
   \hline
  \; $N^{\phantom{|}}$ & $\beta(N)$  \\[.0cm]
   \hline
   1\;& \;${1}\cdot$              \\
   2\; & \;\:\,\,$1{0}\cdot01$     \\
   3\; & \;$10{0}\cdot01$          \\
   4\; & \;$10{1}\cdot01$          \\
   5\; & \;\:\,$100{0}\cdot1001$   \\
   6\; & \;\:\,$101{0}\cdot0001$   \\
   7\; & \;$1000{0}\cdot0001$      \\
   8\; & \;$1000{1}\cdot 0001$     \\
   \hline
 \end{tabular}\quad
 \begin{tabular}{|r|c|}
   \hline
  \; $N^{\phantom{|}}$ & $\beta(N)$  \\[.0cm]
   \hline
   9\; & \;$1001{0}\cdot0101$    \\
   10\; & \;$1010{0}\cdot0101$   \\
   11\; & \;$1010{1}\cdot0101$   \\
   12\; & \;\,\,$10000{0}\cdot101001$      \\
   13\; & \;\,\,$10001{0}\cdot001001$      \\
   14\; & \;\,\,$10010{0}\cdot001001$      \\
   \phantom{a}15\; & \;\,\,$10010{1}\cdot001001$   \\
   16\; & \;\, $10100{0}\cdot100001$\\
   \hline
 \end{tabular}\quad
  \begin{tabular}{|r|c|}
   \hline
   \;$N^{\phantom{|}}$ & $\beta(N)$ \\[.0cm]
   \hline
   17\; & \;\,\,\,$101010\cdot000001$ \\
   18\; & \;$1000000\cdot000001$   \\
   19\; & \;$1000001\cdot000001$   \\
   20\; & \;$1000010\cdot010001$   \\
   21\; & \;$1000100\cdot010001$   \\
   22\; & \;$1000101\cdot010001$   \\
   23\; & \;$1001000\cdot100101$   \\
   \phantom{a}24\; & \;$1001010\cdot000101$  \\
   \hline
 \end{tabular}

\medskip

\section{The recursive structure theorem}\label{sec:3}

The Lucas numbers $(L_n)=(2, 1, 3, 4, 7, 11, 18, 29, 47, 76,123, 199, 322,\dots)$ are defined by
$$  L_0 = 2,\quad L_1 = 1,\quad L_n = L_{n-1} + L_{n-2}\quad {\rm for \:}n\ge 2.$$
The Lucas numbers have a particularly simple base phi representation.

\noindent From  the well-known formula
$L_{2n}=\varphi^{2n}+\varphi^{-2n}$, and the recursion $L_{2n+1}=L_{2n}+L_{2n-1}$ we have for all $n\ge 1$
$$ \beta(L_{2n}) = 10^{2n}\cdot0^{2n-1}1,\quad \beta(L_{2n+1}) = 1(01)^n\cdot(01)^n.$$
By iterated application of the double carry and the golden mean shift to $\beta(L_{2n+1})+\beta(1)$,\; we find that for all $n\ge 1$
$$\beta(L_{2n+1}+1) = 10^{2n+1}\cdot(10)^n01.$$

\noindent As in \cite{Dekk-phi-FQ} we partition the natural numbers into Lucas intervals\: 
$$\Lambda_{2n}:=[L_{2n},\,L_{2n+1}] \quad{\rm and\quad} \Lambda_{2n+1}:=[L_{2n+1}+1,\, L_{2n+2}-1].$$
The basic idea behind this partition is that if
 $$\beta(N) = d_{L}d_{L-1}\dots d_1d_0\cdot d_{-1}d_{-2} \dots d_{R+1}d_R,$$
then the left most index $L=L(N)$ and the right most index $R=R(N)$ satisfy
$$L(N)=|R(N)|=2n \;{\rm iff}\; N\in \Lambda_{2n}, \quad L(N)=2n\!+1.\; |R(N)|=2n\!+2 \;{\rm iff}\; N\in \Lambda_{2n+1}.$$
This is not hard to see from the simple expressions we have for the $\beta$-representations of the Lucas numbers, see also Theorem 1 in \cite{Grabner94}.

\medskip

In some sense, odd Lucas intervals are not small enough. 
 To obtain recursive relations, the interval $\Lambda_{2n+1}=[L_{2n+1}+1, L_{2n+2}-1]$ has to be divided into three subintervals. These three intervals are\\[-.8cm]
 \begin{align*}
I_n:=&[L_{2n+1}+1,\, L_{2n+1}+L_{2n-2}-1],\\
J_n:=&[L_{2n+1}+L_{2n-2},\, L_{2n+1}+L_{2n-1}],\\
K_n:=&[L_{2n+1}+L_{2n-1}+1,\, L_{2n+2}-1].
\end{align*}

\noindent Note that $I_n$ and $K_n$ have the same length $L_{2n-2}-1$, and that $J_n$ has length $L_{2n-3}+1$.

\noindent It will be very convenient to use the free group versions of words of 0's and 1's. So, for example, $(01)^{-1}0001=1^{-1}001$.

\begin{theorem*}{\bf [Recursive structure theorem]}\label{th:rec} 

\noindent{\,\bf I\;} For all $n\ge 1$ and $k=1,\dots,L_{2n-1}$
one has $ \beta(L_{2n}+k) =  \beta(L_{2n})+ \beta(k) = 10\dots0 \,\beta(k)\, 0\dots 01.$
\noindent{\bf II} For all $n\ge 2$ and $k=1,\dots,L_{2n-2}-1$\\[-.8cm]
\begin{align*}
I_n:&\quad \beta(L_{2n+1}+k) = 1000(10)^{-1}\beta(L_{2n-1}+k)(01)^{-1}1001,\\[-.1cm] K_n:&\quad\beta(L_{2n+1}+L_{2n-1}+k)=1010(10)^{-1}\beta(L_{2n-1}+k)(01)^{-1}0001.
\end{align*}\\[-.8cm]
Moreover, for all $n\ge 2$ and $k=0,\dots,L_{2n-3}$\\[-.4cm]
$$\hspace*{0.7cm}J_n:\quad\beta(L_{2n+1}+L_{2n-2}+k) = 10010(10)^{-1}\beta(L_{2n-2}+k)(01)^{-1}001001.$$
\end{theorem*}

\medskip

\noindent {\it Proof.}\\ 
{\bf I} As noted in \cite{Dekk-phi-FQ}, Part {\bf I} follows in a very simple way,
 because adding  $\beta(k)$  to $\beta(L_{2n})$ does not give any double carries or golden mean shifts, when  $k$ is smaller than $L_{2n-1}$. \\
{\bf II Part $I_n$.}\\ Fix a number $k$ with $k\in \{1,\dots,L_{2n-2}-1\}$. Write, with $L=2n-1, R=-2n$,
$$\beta(L_{2n-1}+k)= 10d_{L-2}\dots d_0\cdot d_{-1}\dots d_{R+2} 01.$$
Let us write $\beta(L_{2n}) = 10^{2n}\cdot0^{2n-1}1$ as
$$\beta(L_{2n}) = 10e_{L-1}\dots  e_0\cdot e_{-1}\dots e_{R+2}01,$$
where all $e_i$ are equal to 0.
From $L_{2n+1}+k = L_{2n-1}+k +L_{2n}$ we have $\beta(L_{2n+1}+k) \doteq\beta(L_{2n-1}+k) +\beta(L_{2n})$, and so,
since $d_L+e_L=1,\,d_{L-1}+e_{L-1}=0, \,d_{R+1}+e_{R+1}=0,$ and $\,d_{R}+e_{R}=2$, one obtains
\begin{align*}
\beta(L_{2n+1}+k) \doteq &\; 110\,d_{L-2}\dots d_0\cdot d_{-1}\dots d_{R+2}\,02\\
=&\; 1000\,d_{L-2}\dots d_0\cdot d_{-1}\dots d_{R+2}1001\\
=&\; 1000(10)^{-1}\,\beta(L_{2n-1}+k)(01)^{-1}1001.\\[-.9cm] 
\end{align*} 
{\bf II  Part $J_n$.}\\ Note first that\\[-.3cm] 
$$\beta( 2L_{2n})\doteq 2\,0^{2n}\cdot0^{2n-1}\,2 =1001\,0^{2n-2}\cdot  0^{2n-2}\,1001=100\,\beta(L_{2n-2})\,1^{-1}0\,1001.$$
We now exploit the equation $L_{2n+1} + L_{2n-2}= 2L_{2n}$. Fix a number $k$ with $k\in \{0,\dots,L_{2n-3}\}$. Then\\[-.8cm]
\begin{align*} 
\beta(L_{2n+1} + L_{2n-2}+k)=&\; \beta( 2L_{2n}+k)\doteq\beta( 2L_{2n})+\beta(k)\\
                        =&\; 100\,\beta(L_{2n-2})\,1^{-1}0\,1001+\beta(k)\\
                        =&\;  100\beta(L_{2n-2}+k)1^{-1}01001,
\end{align*}
where we used Part {\bf I} in the second and in the last step
(in the second step a version of Part {\bf I}, with $L_{2n}$ replaced by $2L_{2n}$).
Part $J_n$ is now proved, since    $10010(10)^{-1} =100$, and $(01)^{-1}001001=1^{-1}01001$.

\noindent {\bf II Part $K_n$.}\\
This is more involved than the proofs of Part $I_n$ and Part $J_n$. We first give a one line proof of the following formula, which can also be found  in Lemma 3.3 of \cite{Hart-Sanchis-99}. For all $n\ge 1$
$$\beta(L_{2n}-1)=(10)^n\cdot 0^{2n-1}1.$$
The proof is a matter of applying the golden mean shift $n$ times:
$$\beta(L_{2n}-1)+\beta(1)\doteq(10)^{n-1}11\cdot^{2n-1}1\doteq(10)^{n-2}1100\cdot 0^{2n-1}1=\ldots=1(00)^n\cdot 0^{2n-1}1=\beta(L_{2n}).$$
The result follows then from the unicity of $\beta$-representations.\\
Note that we proved the $K_n$-formula
$$\beta(L_{2n+1}+L_{2n-1}+k)=10\beta(L_{2n-1}+k)(01)^{-1}0001$$
for $k=L_{2n-2}-1$, since
$$L_{2n+1}+L_{2n-1}+L_{2n-2}-1=L_{2n+1}+L_{2n}-1=L_{2n+2}-1, \quad L_{2n-1}+L_{2n-2}-1=L_{2n}-1,$$
and
$$\beta(L_{2n+2}-1)=(10)^{n+1}\cdot 0^{2n+1}1=10\,(10)^n\cdot 0^{2n-1}\,11^{-1}\,001=10\beta(L_{2n}-1)(01)^{-1}0001.$$
Next, note that $k=L_{2n-2}-1$ in the right side of the $K_n$-formula gives the last element $L_{2n-1}+L_{2n-2}-1=L_{2n}-1$ of the Lucas interval $\Lambda_{2n-1}=[L_{2n-1}+1, L_{2n}-1]$, and $k=1$ gives the first element of $\Lambda_{2n-1}$. This implies that repeatedly adding $\beta(1)$ to this first element gives  $\beta$-representations all restricted to the same range. But since we proved the correctness of the $K_n$-formula for the last number obtained, the formula must then also be correct for all previous numbers, again by uniqueness of the $\beta$-representations.  $\Box$

\noindent AMS Classification Numbers: 11D85, 11A63, 11B39

\end{document}